\documentclass[conference,a4paper]{IEEEtran}

\usepackage{verbatim}
\usepackage{amsfonts,amsmath,mathrsfs,amssymb,amsbsy}
\usepackage[final]{graphicx}
\usepackage{times,cite}

\long\def\symbolfootnote[#1]#2{\begingroup%
\def\thefootnote{\fnsymbol{footnote}}\footnote[#1]{#2}\endgroup}

\usepackage{verbatim}
\usepackage[]{float,latexsym}
\usepackage{amsfonts,amsmath,mathrsfs,amssymb,amsbsy}
\usepackage{url}

\newtheorem{theorem}{Theorem}[section]

\newcommand{\Prob}{\mathsf{P}}
\newcommand{\Expect}{\mathsf{E}}

\usepackage{color}
\definecolor{lightblue}{rgb}{.7, .8, 1}
\definecolor{lightgreen}{rgb}{.6, 1, .6}
\usepackage{color}
\definecolor{brown}{rgb}{1,0.38,0.03}

\definecolor{OliveGreen}{rgb}{.2,0.6,0.2}
\definecolor{BrickRed}{rgb}{.7,0.2,0.2}

\newcommand{\ignore}[1]{} 

\long\def\symbolfootnote[#1]#2{\begingroup%
\def\thefootnote{\fnsymbol{footnote}}\footnote[#1]{#2}\endgroup}

\DeclareMathOperator*{\argmax}{arg\,max}

\newcommand{\taug}{\tau_{\scriptscriptstyle \text{G}}}

\begin{document}

\sloppy

\title{Non-parametric Quickest Change Detection for Large Scale Random Matrices}

\author{
  \IEEEauthorblockN{Taposh Banerjee, Hamed Firouzi and Alfred O. Hero III}\\
  \IEEEauthorblockA{Department of EECS\\
    University of Michigan\\
    Ann Arbor, MI, USA\\
    Email: {taposh, firouzi, hero}@umich.edu}
}



\maketitle

\begin{abstract}
The problem of quickest detection of a change in the distribution of a $n\times p$ random matrix based on a sequence of observations having a single unknown change point is considered. The forms of the pre- and post-change distributions of the rows of the matrices are assumed to belong to the family of elliptically contoured densities with sparse dispersion matrices but are otherwise unknown.  We propose a non-parametric stopping rule that is based on a novel summary statistic related to k-nearest neighbor correlation between columns of each observed random matrix. In the large scale regime of $p\rightarrow \infty$ and $n$ fixed we show that, among all functions of the proposed summary statistic,  the proposed stopping rule is asymptotically optimal under 
a minimax quickest change detection (QCD)  model.

\end{abstract}

\section{Introduction} \label{sec:Intro}
In this paper we consider the problem of sequential detection of 
a change in the distribution of a sequence of large scale random matrices. 
The random matrices have i.i.d. rows where the pre-change and post-change 
distributions of the rows are known to belong to the elliptically contoured family but are otherwise unknown. 
This large scale non-parametric sequential detection 
problem has applications in multivariate time-series analysis, stochastic finance, social networks 
and failure detection, among others. In multivariate time-series analysis, it is of interest to know if the coefficients 
of the time series has changed over time. In stochastic finance, it is of interest 
to detect a sudden change in the correlation between a set of stocks being monitored.  
In social networks, it is of interest to detect an abrupt change in the interaction level between a pair of agents. 
In failure detection, often the dynamics of a mechanical structure can be characterized by multi-variate data, 
and a change in the dynamics should be detected as quickly as possible. 

In such cases the observations can be described as a sequence of random matrices. 
The rows of these random matrices may correspond to approximately independent 
realizations of $p$ different variables, e.g., sampled over blocks of time or sampled in a sequence of repeated experiments.  
For example, in the case of detecting a change in the coefficients of a Gaussian univariate time series, 
$p$ successive time samples may be acquired over $n$ well separated blocks of time. A change in the coefficients of the time series is reflected in a change in the correlation matrix associated with each block. 
In stochastic finance, we may have access to multiple instances of stock values over a day or week, 
and a change in correlation may occur only at the end of the day or week. 

In this paper we consider the problem of quickest detection of a change in 
population dispersion (or correlation) matrix under the assumption of elliptically contoured distribution of the rows of the sequence of $n\times p$ random matrices. 
The results in this paper hold for the big data regime of $p \gg n$ 
 for which $p \rightarrow \infty$ and $n$ is fixed and small. 
The precise mathematical problem is stated in Section~\ref{sec:Prob}. 

If a parametric model for the data is known before and after change, then various efficient procedures 
from the quickest change detection literature (see, e.g., \cite{veer-bane-elsevierbook-2013}, \cite{poor-hadj-qcd-book-2009}, 
and \cite{tart-niki-bass-2014}) can be used for detection. 
However, in the absence of a parametric model, a situation common in Big Data settings,  
no optimal procedures are known. In this paper we propose 
a technique for quickest change detection in this setting. 

Specifically, we propose a novel summary statistic
for the data matrix: the minimal $k$-nearest neighborhood of 
the columns of the random matrix under a correlation magnitude distance.  
We obtain an approximate distribution for the summary statistic in 
the big data regime. 
We show that the distribution of the summary statistic belongs to a one-parameter exponential family, 
with the unknown parameter a function of the underlying distribution of the data matrix. 
 We then 
treat the sequence of summary statistics as our observation sequence, and apply Lorden's test \cite{lord-amstat-1971}. 
This work is motivated by the theory of correlation screening and correlation mining   
\cite{hero-bala-jasta-2011}, and specifically 
the theory of hub discovery in large scale correlation graphs from \cite{hero-bala-IT-2012}.

\section{Problem Description}\label{sec:Prob}
A decision-maker sequentially acquires samples from a family of distributions of 
$n \times p$ random matrices over time, indexed by $m$, leading to the random matrix sequence $\{\mathbb{X}(m)\}_{m \geq 1}$, 
called data matrices. 
For each $m$ the random matrix $\mathbb{X}(m)$ has the following properties. Each of its 
$n$ rows is an independent identically distributed (i.i.d.) sample 
of a $p$-variate random vector $\mathbf{X}(m)$ $=[X_1(m), \cdots, X_p(m)]^T$ with $p \times 1$ 
mean $\boldsymbol{\mu}_m$ and 
$p\times p$ positive definite dispersion matrix $\mathbf{\Sigma}_m$. 
The random vector $\mathbf{X}(m)$ has an elliptically contoured density, 
also called an elliptical density \cite{ander-mutlstat-1996},
\[
f_{\mathbf{X}(m)}(\mathbf{x}) = g_m((\mathbf{x}-\boldsymbol{\mu}_m)^T\boldsymbol{\Sigma}_m^{-1}(\mathbf{x}-\boldsymbol{\mu}_m)),
\]
for some nonnegative strictly decreasing function $g_m$ on $\mathbb{R}^+$. If $\boldsymbol{\mu}_m=0$
and $\mathbf{\Sigma}_m = I_p$, where $I_p$ is the $p \times p$ identity matrix, then 
the random vector $\mathbf{X}(m)$ is said to have a spherical density. 

The samples $\{\mathbb X(m)\}$ are assumed to be statistically independent. 
For some time parameter $\gamma$ the samples are assumed to have common dispersion 
parameter $ \mathbf{\Sigma}_0$ and function $g_0$ for $m < \gamma$  
and common dispersion parameter $\mathbf{\Sigma}_1 \neq \mathbf{\Sigma}_0$ and function $g_1$ for $m\geq \gamma$. 
$\gamma$ is called the change point and the pre-change and post-change distributions of $\mathbb X(m)$ 
are denoted $f_{\mathbf{X}}^0$ and  $f_{\mathbf{X}}^1$, respectively. 
No assumptions are made about the mean parameter $\boldsymbol{\mu}_m$, and can take different 
values for different $m$. 
More specifically, as the rows  of $\mathbb X(m)$ are i.i.d. realizations of 
the elliptically distributed random variable $\mathbf{X}(m)$, this change-point model is described by: 
\begin{equation}\label{eq:PreChg_PstChg_Den}
\begin{split}
\mathbf{X}(m) &\sim f_{\mathbf{X}}^{0}(\mathbf{x}) = g_0((\mathbf{x}-\boldsymbol{\mu}_m)^T\boldsymbol{\Sigma_0}^{-1}(\mathbf{x}-\boldsymbol{\mu}_m)), \; m < \gamma \\ 
       &\sim f_{\mathbf{X}}^{1}(\mathbf{x}) = g_1((\mathbf{x}-\boldsymbol{\mu}_m)^T\boldsymbol{\Sigma_1}^{-1}(\mathbf{x}-\boldsymbol{\mu}_m)), \; m \geq \gamma.
\end{split}
\end{equation}
At each time point $m$ the decision-maker decides to either stop sampling, 
declaring that the change has occurred, i.e., $m\geq \gamma$, or to continue sampling. 
The decision to stop at time $m$ is only a function of $(\mathbb{X}(1), \cdots, \mathbb{X}(m))$. 
Thus, the time at which the decision-maker decides to stop sampling is a stopping time 
for the matrix sequence $\{\mathbb{X}(m)\}$.
The decision-maker's objective is to detect this change in distribution of the data matrices as quickly as possible, 
subject to a constraint on the false alarm rate.  

The above detection problem is an example of the quickest change detection (QCD) problem.
See \cite{poor-hadj-qcd-book-2009}, \cite{veer-bane-elsevierbook-2013}, and \cite{tart-niki-bass-2014} for 
an overview of the QCD literature. In the QCD problem the objective is to find a stopping time $\tau$ 
on the sequence of data matrices $\{\mathbb{X}(m)\}$, so as to minimize a suitable metric on the delay $(\tau-\gamma)$, 
subject to a constraint on a suitable metric on the event of false alarm $\{\tau < \gamma\}$. 
This paper follows the QCD formulation of Pollak \cite{poll-astat-1985}: 
\begin{equation}\label{prob:Pollak}
\begin{split}
 \min_\tau       \sup_{\gamma\geq 1} & \quad \Expect_\gamma[\tau-\gamma| \tau \geq \gamma] \\
 \mbox{subj. to}  &  \quad \Expect_\infty[\tau] \geq \beta, 
\end{split}
\end{equation}
where $\Expect_\gamma$ is the expectation with respect to the probability measure under which the change occurs at $\gamma$, 
$\Expect_\infty$ is the corresponding expectation when the change never occurs, and $\beta \geq 1$ is a 
user-specified constraint on the mean time to false alarm.

If the pre- and post-change densities $f_{\mathbf{X}}^{0}$ and $f_{\mathbf{X}}^{1}$ are known to the decision maker, 
and $\boldsymbol{\mu}_m$ is constant before and after change, 
then algorithms like the Cumulative Sum (CuSum) algorithm \cite{page-biometrica-1954}, \cite{lord-amstat-1971}, 
\cite{mous-astat-1986}, or the Shiryaev-Roberts (SR) family of algorithms \cite{robe-technometrics-1966}, \cite{poll-astat-1985}, 
\cite{tart-etal-thirdord-2012}, can be used for efficient change detection. Both the CuSum algorithm 
and the SR family of algorithms have strong optimality properties with respect to both the popular 
formulations of Lorden \cite{lord-amstat-1971} and that of Pollak \cite{poll-astat-1985}, used in this paper. 

If only the pre-change and post-change functions $g_0$ and $g_1$ are known 
then \eqref{prob:Pollak} is a parametric QCD problem. 
In this case, under the assumption that $\boldsymbol{\mu}_m=\boldsymbol{\mu}_0$, $m < \gamma$, and $\mathbf{\Sigma}_0$ are known, 
efficient QCD algorithms can be designed, having strong asymptotic optimality properties,
based on, e.g., the generalized likelihood ratio (GLR) technique \cite{tart-niki-bass-2014}, 
the mixture based technique  \cite{tart-niki-bass-2014}, or 
the nonanticipating estimation based technique \cite{lord-poll-nonantiest-2005}.   
   
In many situations, however, even the pre- and post-change functions $g_0$ and $g_1$ may be unknown. 
This is the non-parametric QCD setting considered in this paper.
While one can use non-parametric QCD tests based on signs and ranks \cite{gord-Poll-signrank-1994}, 
or based on empirical distribution estimates \cite{li-siri-veer-isit-2014}, 
there are no known optimal solutions to \eqref{prob:Pollak} in the non-parametric setting.   

In this paper we provide an asymptotically optimal solution to the minimax QCD problem  
\eqref{prob:Pollak} in the random matrix setting \eqref{eq:PreChg_PstChg_Den} 
using recently developed large scale random matrix theory \cite{hero-bala-IT-2012}.
The solution is optimal in the following sense.  
The theory from \cite{hero-bala-IT-2012} establishes that a certain 
summary statistic, denoted by $V(\mathbb{X})$, derived from an $n\times p$ random matrix $\mathbb{X}$
has a limiting distribution as $p\rightarrow \infty$ for fixed $n$, the so-called "purely high dimensional regime" 
\cite{hero-bala-submitted-2015}. 
This summary statistic is related to the empirical distribution of the vertex degree 
of the correlation graph associated with the thresholded sample correlation matrix.   
Below we show that the distribution of the statistic $V(\mathbb{X})$ converges to a 
parametric distribution in the exponential family in this purely high dimensional regime. 
We then apply the GLR based Lorden's test \cite{lord-amstat-1971} to the sequence of summary statistics
$\{V(\mathbb{X}(m))\}$ to detect the change efficiently.  
Thus, the proposed stopping rule is asymptotically optimal under the Lorden minimax quickest change detection (QCD) model \cite{lord-amstat-1971}, 
and hence also in terms of solving \eqref{prob:Pollak}, among all rules that 
are stopping rules for the proposed summary statistics sequence.

\section{Summary Statistic for the Data Matrix}\label{sec:SumStat}
In this section we define a summary statistic $V(\mathbb{X})$ and 
then use the results from \cite{hero-bala-IT-2012} to obtain its asymptotic density in the 
purely high dimensional regime of  $p\rightarrow \infty$, $n$ fixed. This asymptotic 
distribution is  a member of 
a one-parameter exponential family. 

For an elliptically distributed random data matrix $\mathbb{X}$ we write
$$\mathbb{X} = [\mathbf{X}_1, \cdots, \mathbf{X}_p] = [\mathbf{X}^T_{(1)}, \cdots, \mathbf{X}^T_{(n)}]^T,$$
where $\mathbf{X}_i = [X_{1i}, \cdots, X_{ni}]^T$ is the $i^{th}$ column and $\mathbf{X}_{(i)} = [X_{i1}, \cdots, X_{ip}]$  
is the $i^{th}$ row. 
Define the sample covariance matrix as 
$$ \mathbf{S} = \frac{1}{n-1} \sum_{i=1}^n (\mathbf{X}_{(i)} - \bar{\mathbf{X}})^T (\mathbf{X}_{(i)} - \bar{\mathbf{X}}),$$
where $\bar{\mathbf{X}}$ is the sample mean of the $n$ rows of $\mathbb{X}$.
Also define the sample correlation matrix as
$$ \mathbf{R} = \mathbf{D_S}^{-1/2} \mathbf{S} \mathbf{D_S}^{-1/2},  $$
where $\mathbf{D_A}$ denotes the matrix obtained by zeroing out all but the diagonal elements of the matrix $\mathbf{A}$.  
Note that, under our assumption that the dispersion matrix $\mathbf{\Sigma}$ of the rows of $\mathbb X$ is positive definite, $\mathbf{D_S}$ is invertible with probability one.
Thus $\mathbf{R}_{ij}$, the element in the $i^{th}$ row and the $j^{th}$ column of the matrix $\mathbf{R}$, is the 
sample correlation coefficient between the $i^{th}$ and $j^{th}$ columns of $\mathbb X$.   

Define $d^{(k)}_{\text{NN}}(i)$ to be the 
sample correlation between the $i$-th column of $\mathbb X$ and 
its $k$-th nearest neighbor in the columns of $\mathbb X$ (in terms of Euclidean distance):
$$d^{(k)}_{\text{NN}}(i) := k^{th} \mbox{ largest order statistic of } \{ |\mathbf{R}_{ij}|; j \neq i\}.$$
Then for fixed $k$, define the summary statistic
\begin{equation}\label{eq:Sum_Stat_V}
V_k(\mathbb{X}):= \max_i  d^{(k)}_{\text{NN}}(i).
\end{equation}
Below we show that the distribution of the statistic $V_k$ can be related to the distribution 
of an integer valued random variable $N_{\delta, \rho}$ which we define below.

For a threshold parameter $\rho \in [0,1]$ define the correlation graph $\mathcal{G}_\rho(\mathbf{R})$ 
associated with the correlation matrix $\mathbf{R}$ 
as an undirected graph with $p$ vertices, each representing a column of the data matrix $\mathbb{X}$. 
An edge is present between vertices $i$ and $j$ if the magnitude of the sample correlation coefficient 
between the $i^{th}$ and $j^{th}$ components of the random vector $\mathbf{X}$ 
is greater than $\rho$, i.e., if $|\mathbf{R}_{ij}| \geq \rho$, $i \neq j$. 
We define $\delta_i$ to be the degree of vertex $i$ in the graph $\mathcal{G}_\rho(\mathbf{R})$. 
For a positive integer $\delta \leq p-1$ we say that a vertex $i$ in the graph $\mathcal{G}_\rho(\mathbf{R})$
is a hub of degree $\delta$ if $\delta_i \geq \delta$. We denote 
by $N_{\delta, \rho}$ the total number of hubs in the correlation graph $\mathcal{G}_\rho(\mathbf{R})$, i.e., 
$$N_{\delta, \rho} = \text{card} \{i: \delta_i \geq \delta \}. $$ 
 
The events $\{V_\delta(\mathbb{X}) \geq \rho\}$ and $\{N_{\delta, \rho} > 0\}$ are equivalent. Hence
\begin{equation}\label{eq:V_eq_N}
\Prob(V_\delta(\mathbb{X}) \geq \rho) = \Prob(N_{\delta, \rho} > 0). 
\end{equation}  

An asymptotic approximation to the probability $\Prob(N_{\delta, \rho} > 0)$ is obtained in \cite{hero-bala-IT-2012} by relating 
$N_{\delta, \rho}$ to a Poisson random variable in the purely high dimensional limit as $p \rightarrow \infty$ and $n $ fixed.
We summarize the approximation in the theorem below.  
We say that a matrix is row sparse of degree $k$ if there are no more than $k$ nonzero entries in any row. We say that 
a matrix is block sparse of degree $k$ if the matrix can be reduced to block diagonal form having a single $k \times k$ block, 
via row-column permutations. 
\medskip
\begin{theorem}[\cite{hero-bala-IT-2012}]\label{thm:CorrScr}
Let $\mathbf{\Sigma}$ be row sparse of degree $k=o(p)$. Also let $p \to \infty$ and $\rho = \rho_p \to 1$ such that 
$p^{1/\delta} (p-1)(1-\rho^2)^{(n-2)/2} \to e_{n,\delta} \in (0, \infty)$. 
\begin{enumerate}
\item $$\Prob(N_{\delta, \rho} > 0) \to 1-\exp(-\Lambda J_\mathbf{{X}}/\phi(\delta)),$$
where 
$$\Lambda = \lim_{p\to \infty, \rho \to 1} \Lambda(\rho) = ((e_{n,\delta} a_n ) / (n-2) )^{\delta}/\delta!,$$
with
$$\Lambda(\rho) = p {p-1 \choose \delta} P_0(\rho)^\delta,$$
$$P_0(\rho) = a_n \int_{\rho}^1 (1-u^2)^{\frac{n-4}{2}} du,$$
$$a_n = 2B((n-2)/2, 1/2) \mbox{ with } B(l,m) \mbox{ the beta function}, $$
$\phi(\delta)=2$ if $\delta = 1$, $\phi(\delta)=1$ otherwise, 
and $J_\mathbf{X}$ is a positive real number that is a function of the joint density of $\mathbf{X}$. 
\item If the dispersion matrix $\mathbf{\Sigma}$ of the p-variate vector $\mathbf{X}$ is block sparse of degree $k$, 
then 
$$J_\mathbf{{X}}=1 + O((k/p)^{\delta+1}).$$ 
In particular, if the dispersion matrix $\mathbf{\Sigma}$ is diagonal then $J_\mathbf{{X}}=1$.
\end{enumerate} 
\end{theorem}
\medskip

Using \eqref{eq:V_eq_N} and Theorem~\ref{thm:CorrScr}, 
the  large $p$  distribution of $V_k$ defined in \eqref{eq:Sum_Stat_V} can be approximated, for $k=\delta$, by
\begin{equation}\label{eq:CDF_V}
\Prob(V_\delta(\mathbb{X}) \leq \rho) = \exp(-\Lambda(\rho) J_\mathbf{{X}}/\phi(\delta)), \; \rho \in [0,1], 
\end{equation}
where $\Lambda(\rho)$ is as defined in Theorem~\ref{thm:CorrScr}.  
Although the theorem is valid for large values of $\rho$, numerical experiments \cite{hero-bala-IT-2012}
have shown that the approximation remains accurate for smaller values of $\rho$ as long as $n$ is small and $p\gg n$. 

The distribution \eqref{eq:CDF_V} is differentiable everywhere except at $\rho=0$ 
since $P(V_\delta(\mathbb{X})=0)>0$ when using the finite $p$ and $\rho<1$ approximation $\Lambda_\rho$ for $\Lambda$ specified in 
Theorem~\ref{thm:CorrScr}. 
For $\rho>0$ and large $p$,  $V_\delta$ has density  
\begin{equation}\label{eq:PDF_V_1}
f_V(\rho) = - \frac{\Lambda'(\rho) }{\phi(\delta)}J_\mathbf{{X}} \exp\left(-\frac{\Lambda(\rho)}{\phi(\delta)} J_\mathbf{{X}}\right), \; \rho \in (0,1]. 
\end{equation}
Note that $f_V$ in \eqref{eq:PDF_V_1} is the density of the Lebesgue continuous component of 
the distribution \eqref{eq:CDF_V} and that it integrates to $1-O(e^{-p^2})$ over $\rho\in (0,1]$.

The density $f_V$ is a member of a one-parameter exponential family with $J_\mathbf{X}$ as the unknown parameter.  
This follows from the relations below. First 
\begin{equation}\label{eq:Lambdarho}
\begin{split}
\Lambda(\rho) &=  p {p-1 \choose \delta} \left(a_n \int_{\rho}^1 (1-u^2)^{\frac{n-4}{2}} du\right)^\delta \\
&= C \; T(\rho)^\delta,
\end{split}
\end{equation}
where 
\begin{equation}\label{eq:C_pndelta}
C=C_{p,n,\delta}=p {p-1 \choose \delta} a^\delta_n
\end{equation}
does not depend on $\rho$, and 
\begin{equation}\label{eq:Trho}
T(\rho) = \int_{\rho}^1 (1-u^2)^{\frac{n-4}{2}} du.
\end{equation}
Using \eqref{eq:Lambdarho} and noting that $T(\rho)' = -(1-\rho^2)^{\frac{n-4}{2}}$, we have for $ \rho \in [0,1]$,
the exponential family form of the density $f_V$ with parameter $J_\mathbf{{X}}$:
\begin{equation}\label{eq:PDF_V}
\begin{split}
f_V &(\rho; J_\mathbf{{X}})  \\
=& \frac{C \delta}{\phi(\delta)} T(\rho)^{\delta-1} (1-\rho^2)^{\frac{n-4}{2}} J_\mathbf{{X}} \exp\left(-\frac{C  T(\rho)^\delta}{\phi(\delta)} J_\mathbf{{X}}\right).
\end{split}
\end{equation}

The constant $\delta$ in \eqref{eq:PDF_V} is a fixed design parameter 
that can be selected to maximize change detection performance according to \eqref{prob:Pollak}. 
In the sequel, we fix $\delta=1$. 
For this value of $\delta$, the statistic $V_\delta$ reduces to the nearest neighbor (correlation) distance
\begin{equation}\label{eq:Sum_Stat_V_deltaeq1}
V(\mathbb{X})= \max_{i \neq j} |\mathbf{R}_{ij}|,
\end{equation}
and the density in \eqref{eq:PDF_V} reduces to
\begin{equation}\label{eq:PDF_V_delta1}
f_V(\rho; J) = \frac{C}{2} (1-\rho^2)^{\frac{n-4}{2}} J \exp\left(-\frac{C}{2} J \;T(\rho) \right), \; \rho \in (0,1],
\end{equation}
where we have suppressed subscript $\mathbf{X}$ in the exponential family parameter $J$ on the distribution of $\mathbf{X}$. 

\begin{figure}[htb]
\center \vspace{-0.3cm}
\includegraphics[width=8cm, height=5cm]{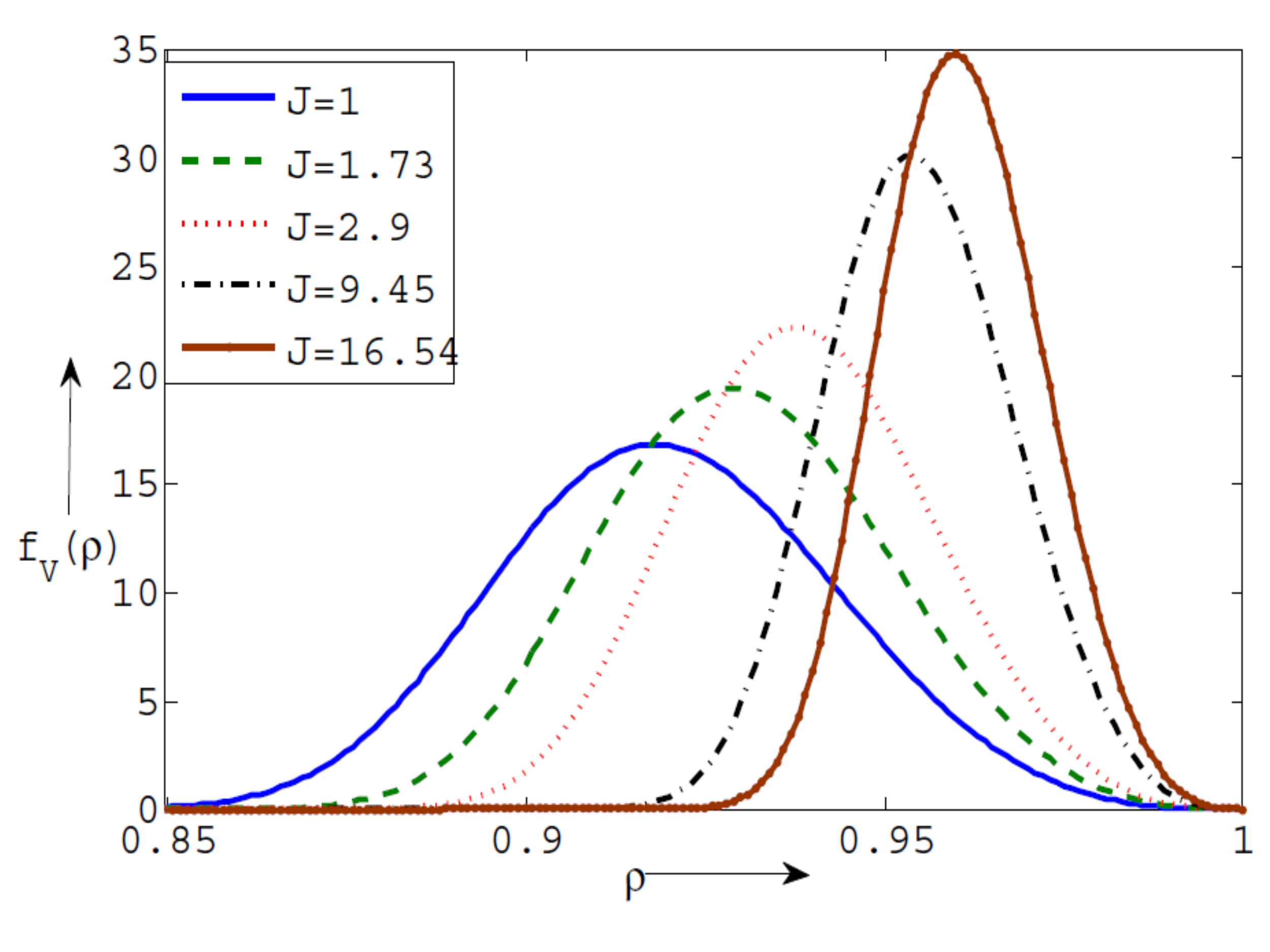}\vspace{-0.5cm}
\caption{Plot of density $f_V$ in \eqref{eq:PDF_V} for various values of the parameter $J$ for $n=10, p=100$. 
 This is the density of the summary statistic used to detect the change in covariance of the random matrix sequence $\mathbb X$.}
\label{fig:fV_density}\vspace{-0.3cm}
\end{figure}

In Fig.~\ref{fig:fV_density} is plotted the density $f_V$ for various values of $J$ for 
$n=10$, and $p=100$. We note that for the chosen values of $n$ and $p$, the density is concentrated close to $1$,
consistent with large values of $\rho$ arising in the purely high dimensional regime assumed in Theorem~\ref{thm:CorrScr}. 

\section{QCD for large scale random matrices}
Here we apply the asymptotic results derived in Section~\ref{sec:SumStat} to quickest 
change detection of the distribution of the summary statistic $V$. 
Assume that both the pre- and post-change dispersion matrices, 
$\mathbf{\Sigma}_0$ and $\mathbf{\Sigma}_1$, are row sparse with degree $k=o(p)$, 
and map the data matrix sequence to the sequence of summary statistics $\{V_\delta(\mathbb{X}(m))\}_{m \geq 1}$. 
For simplicity we refer to this sequence by $\{V(m)\}$. 
Let $J_0$ and $J_1$ be the value of parameter $J$ before and after change, respectively.  
The QCD problem on the density $f_{\mathbf{X}}$, depicted in \eqref{eq:PreChg_PstChg_Den}, 
is reduced to the  QCD problem on the density $f_V$:   
\begin{equation}
\begin{split}
V(m) &\sim f_V(\cdot; J_0), \; m < \gamma \\ 
    &\sim f_V(\cdot; J_1), \; m \geq \gamma.
\end{split}
\end{equation}

We recall from Theorem~\ref{thm:CorrScr} that if the dispersion matrix $\mathbf{\Sigma}_0$ is diagonal 
then $J_0=1$. Thus, if the pre-change dispersion matrix is diagonal, then the QCD problem reduces to the 
parametric QCD problem with unknown post-change parameter $J$: 
\begin{equation}\label{eq:QCDProb_1toother}
\begin{split}
V(m) &\sim f_V(\cdot; 1), \quad\hspace{1.25cm} m < \gamma \\ 
    &\sim f_V(\cdot; J), \; \quad J \neq 1, \; m \geq \gamma.
\end{split}
\end{equation}
If the dispersion matrix $\mathbf{\Sigma}_0$ is only block sparse with degree $k \ll p$, 
by assertion 2 of Theorem~\ref{thm:CorrScr}, we can use the approximation $J_0 \approx 1$. 

Consider the following QCD test, defined by the stopping time $\taug$:
\begin{equation}\label{eq:taug}
\begin{split}
\taug =\inf_{m \geq 1}  \left\{\max_{1 \leq \ell \leq m} \sup_{J: |J-1| \geq \epsilon} \sum_{i=\ell}^m \log \frac{f_V(V(i);J)}{f_V(V(i); 1)} > A\right\},
\end{split}
\end{equation}
where $A$ and $\epsilon>0$ are user-defined parameters. 
The parameter $A$ is a threshold used to control the false alarm rate.
The parameter $\epsilon$ represents the minimum magnitude of change, away from $J=1$, that the user wishes to detect. 

The stopping rule $\taug$ was shown to be asymptotically optimal in \cite{lord-amstat-1971} for a related QCD problem
when 1) the marginal density $f_V(v;\cdot)$ of the observation sequence $\{V(m)\}$ 
is of known form that is a member of a one-parameter exponential family and 2) when the parameter $J_0$  of the pre-change density is known. 
Both of these properties are satisfied for the summary statistic $V=V(\mathbb X)$ for 
the QCD model in \eqref{eq:QCDProb_1toother} defined above, since $J_0=1$. 
Due to the results in \cite{lai-ieeetit-1998}, the stopping rule $\taug$ is asymptotically 
optimal for the problem in \eqref{prob:Pollak} as well.  

The following theorem establishes strong asymptotic optimality of this test. 
\medskip
\begin{theorem}[\cite{lord-amstat-1971}, \cite{lai-ieeetit-1998}]\label{thm:Lorden}
Fix any $\epsilon > 0$.
\begin{enumerate}  
\item For the stopping rule $\taug$, the supremum in \eqref{prob:Pollak} is 
achieved at $\gamma=1$, i.e., 
\[
\sup_{\gamma\geq 1} \; \Expect_\gamma[\taug-\gamma| \taug \geq \gamma] = \Expect_1[\taug-1].
\]
\item Setting $A=\log \beta$ ensures that as $\beta \to \infty$,  
\[
\Expect_\infty[\taug] \geq \beta (1+o(1)),
\]
and for each possible true post-change parameter $J$, with $|J-1|\geq \epsilon$,
\begin{equation}
\begin{split}
\Expect_1[\taug] &= \frac{\log \beta}{I(J)} (1+o(1)) \\
&= \inf_{\tau: \Expect_\infty[\tau] \geq \beta}\sup_{\gamma\geq 1} \; \Expect_\gamma[\tau-\gamma| \tau \geq \gamma] (1+o(1)),
\end{split}
\end{equation}
where $I(J)$ is the Kullback-Leibler divergence between the densities $f_V(\cdot; J)$ and $f_V(\cdot; 1)$. 
\end{enumerate}
\end{theorem}
\medskip
Theorem~\ref{thm:Lorden} implies that the stopping rule $\taug$ is uniformly asymptotically optimal for each post-change parameter $J$, 
as long as $|J-1| \geq \epsilon$. For convenience of implementation one can also use the window limited variation 
of $\taug$ as suggested in \cite{lai-ieeetit-1998}. 

\vspace{-0.3cm}
\section{Numerical Results}
Here we apply the stopping rule $\taug$ in \eqref{eq:taug} to the problem of detecting a change 
in the distribution when the $\{\mathbb X(m)\}$ are Gaussian distributed random matrices. 
In this case the dispersion $\mathbf{\Sigma}$ is the covariance matrix of the rows of $\mathbb X$. 
The pre-change covariance is the $p\times p$ diagonal matrix $\mathbf{\Sigma}_0=\text{diag}(\sigma^2_i)$, where $\sigma_i^2 >0$ are arbitrary component-wise variances. 
The post-change covariance matrix $\mathbf{\Sigma}_1$ is obtained by replacing the $k \times k$ top left block 
of the identify matrix $\mathbf{I}_p$ by a sample from the Wishart distribution. 
We set $n=10$, $p=100$, and $k=5$. 

To implement $\taug$ we have chosen $\epsilon=1.5$, 
and we use the the maximum likelihood estimator which, 
as a function of $m$ samples $(V(1), \cdots, V(m))$ from $f_V(\cdot, J)$, is given by
\begin{equation}\label{eq:MLEst_J}
\hat{J}(V(1), \cdots, V(m)) = \frac{1}{\frac{C}{2}\frac{1}{m} \sum_{i=1}^m T(V(i))}. 
\end{equation}
Specifically, 
\begin{equation}
\begin{split}
\argmax_{J: J \geq 2.5} \;  \log \sum_{i=\ell}^m & \frac{f_V(V(i);J)}{f_V(V(i); 1)} \\
&= \max\{2.5, \hat{J}(V(\ell), \cdots, V(m))\}.
\end{split}
\end{equation}

\begin{figure}[htb]
\center
\includegraphics[width=8cm, height=5cm]{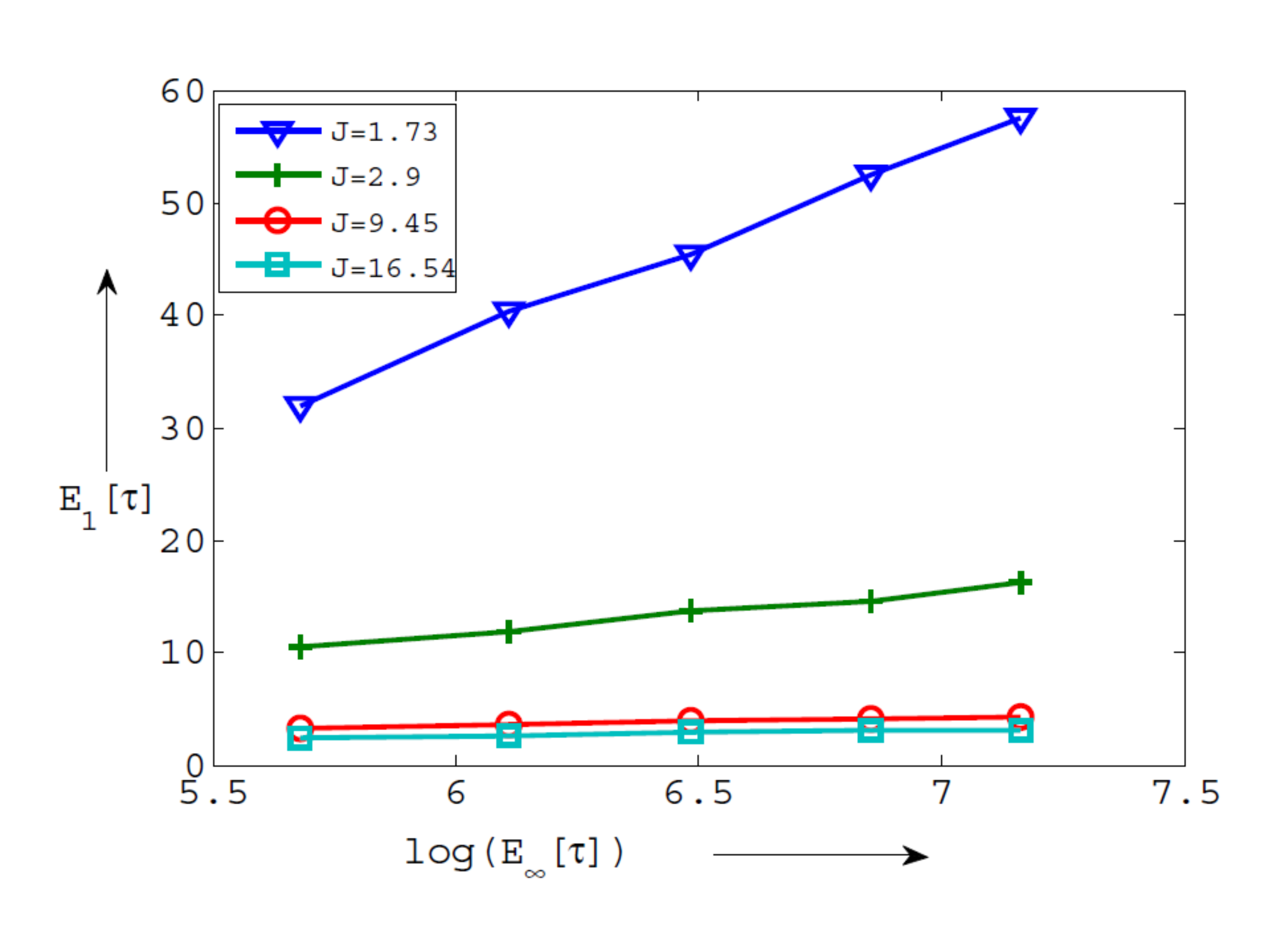}\vspace{-0.2cm}
\caption{The empirical mean time to detect vs  mean time to false alarm (in log scale). 
The mean time to detect decreases as the parameter $J$ increases. }
\label{fig:Delay_MFA}
\end{figure}

In Fig.~\ref{fig:Delay_MFA} we plot the delay ($\Expect_1[\tau]$) vs the log of mean time to false alarm ($\log \Expect_\infty[\tau]$) 
for various values of the post-change parameter $J$. The values in the figure are obtained 
by choosing different values of the threshold $A$ and estimating the delay by choosing the change point $\gamma=1$ 
and simulating the test for $500$ sample paths. 
The mean time to false alarm values are estimated by simulating the test for $1500$ sample paths. 
The parameter $J$ for the post-change distribution is estimated 
using the maximum likelihood estimator \eqref{eq:MLEst_J}.   
As predicted by the theory, the delay vs log of false alarm trade-off curve is approximately linear. 
For larger values of $J$, the Kullback-Leibler (K-L) divergence between $f_V(\cdot, J)$ and $f_V(\cdot, 1)$ is larger, resulting 
in smaller delays. For the chosen values of the post-change parameters $J=1.73$, $2.9$, $9.45$ and $16.54$, 
the corresponding K-L divergence values $I(J)$ are $0.127$, $0.41$, $1.35$ and $1.86$, respectively. 

\begin{figure}[htb]
\center\vspace{-0.3cm}
\includegraphics[width=8cm, height=5cm]{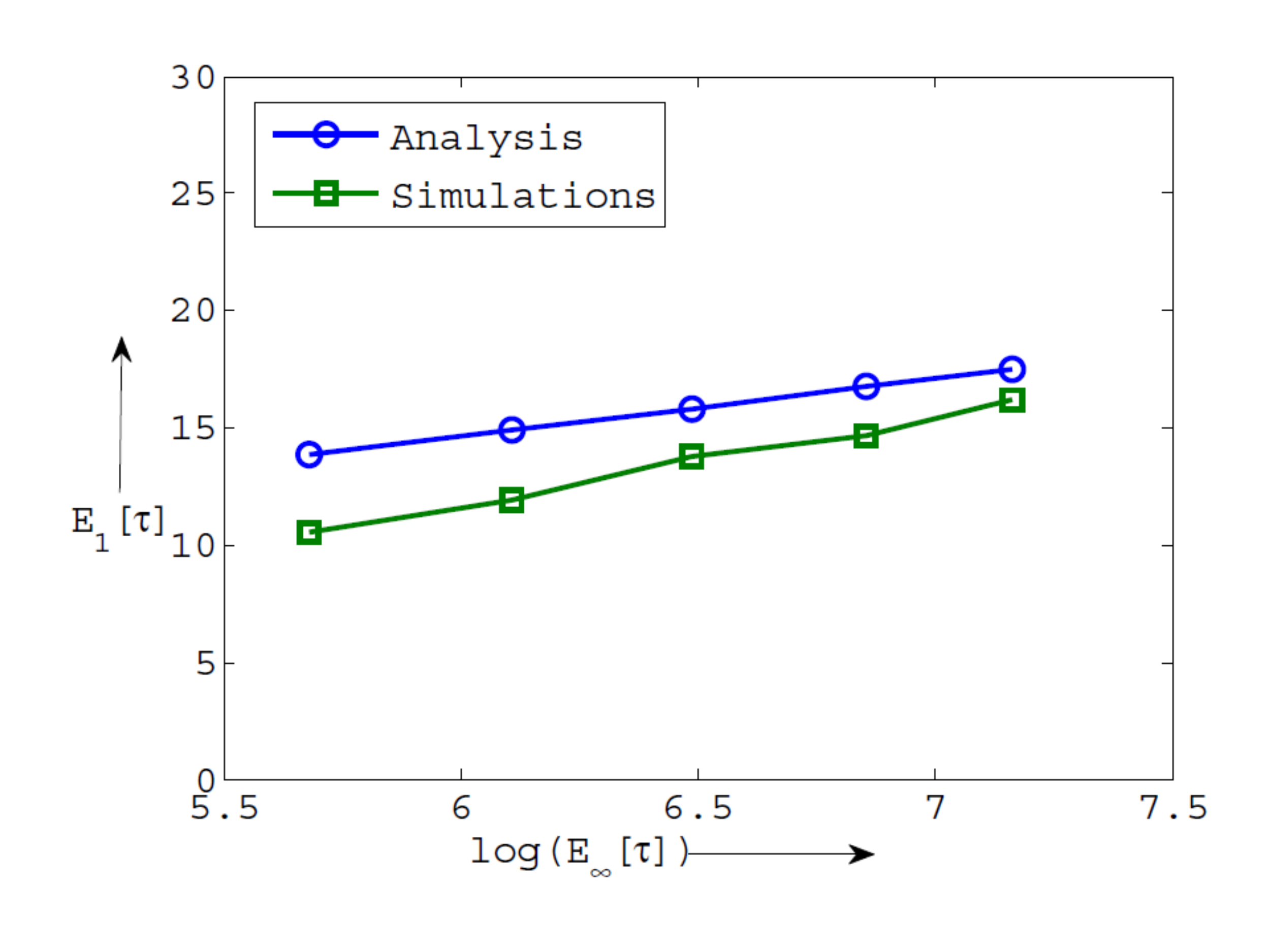}\vspace{-0.3cm}
\caption{Comparison of the delay vs false alarm trade-off curve for $J=2.9$ from Fig.\ref{fig:Delay_MFA} with
the values predicted by the theory: $\frac{\log \Expect_\infty[\tau]}{I(J)} = \frac{\log \Expect_\infty[\tau]}{0.41}$.}
\label{fig:Delay_MFA_AnaSim}\vspace{-0.5cm}
\end{figure}

In Fig.~\ref{fig:Delay_MFA_AnaSim} we compare the delay vs false alarm trade-off curve
for the post-change parameter $J=2.9$ plotted in Fig.~\ref{fig:Delay_MFA}, with the values 
predicted by the theory: $\frac{\log \Expect_\infty[\tau]}{I(J)}$. We see from Fig.~\ref{fig:Delay_MFA_AnaSim}
that the predictions are quite accurate. 
We have obtained similar results when the test was simulated for different block sizes $k$. 
Thus, the change can be efficiently detected using our proposed methodology. 

\section{Conclusions and Future Work}
We have introduced a novel summary statistic based on correlation mining and hub 
discovery for performing non-parametric quickest change detection (QCD) 
on a sequence of large scale random matrices.   The proposed QCD algorithm 
is strongly optimal in the sense of Lorden \cite{lord-amstat-1971} and Pollak \cite{poll-astat-1985}
among all detection algorithms that use our summary statistic.  
Future work will  include extensions to local summary statistics and experiments 
with QCD in  real applications that yield sequences of large scale random matrix measurements.
\vspace{-0.1cm}
\section{Acknowledgments}\vspace{-0.1cm}
This work was partially supported by the Consortium for Verification Technology under Department of Energy National Nuclear Security Administration award number DOE-NA0002534.
\vspace{-0.1cm}


\vspace{-0.3cm}
\bibliographystyle{ieeetr}



\bibliography{QCD_verSubmitted}

\begin{thebibliography}{10}

\bibitem{veer-bane-elsevierbook-2013}
V.~V. Veeravalli and T.~Banerjee, {\em Quickest Change Detection}.
\newblock Elsevier: E-reference Signal Processing, 2013.
\newblock \url{http://arxiv.org/abs/1210.5552}.

\bibitem{poor-hadj-qcd-book-2009}
H.~V. Poor and O.~Hadjiliadis, {\em Quickest detection}.
\newblock Cambridge University Press, 2009.

\bibitem{tart-niki-bass-2014}
A.~G. Tartakovsky, I.~V. Nikiforov, and M.~Basseville, {\em Sequential
  Analysis: {Hypothesis} Testing and Change-Point Detection}.
\newblock Statistics, CRC Press, 2014.

\bibitem{lord-amstat-1971}
G.~Lorden, ``Procedures for reacting to a change in distribution,'' {\em Ann.
  Math. Statist.}, vol.~42, pp.~1897--1908, Dec. 1971.

\bibitem{hero-bala-jasta-2011}
A.~Hero and B.~Rajaratnam, ``Large-scale correlation screening,'' {\em J. Amer.
  Statist. Assoc.}, vol.~106, no.~496, pp.~1540--1552, 2011.

\bibitem{hero-bala-IT-2012}
A.~Hero and B.~Rajaratnam, ``Hub discovery in partial correlation graphs,''
  {\em IEEE Trans. Inf. Theory}, vol.~58, no.~9, pp.~6064--6078, 2012.

\bibitem{ander-mutlstat-1996}
T.~W. Anderson, {\em An Introduction to Multivariate Statistical Analysis}.
\newblock New York, NY: Wiley, 2003.

\bibitem{poll-astat-1985}
M.~Pollak, ``Optimal detection of a change in distribution,'' {\em Ann.
  Statist.}, vol.~13, pp.~206--227, Mar. 1985.

\bibitem{page-biometrica-1954}
E.~S. Page, ``Continuous inspection schemes,'' {\em Biometrika}, vol.~41,
  pp.~100--115, June 1954.

\bibitem{mous-astat-1986}
G.~V. Moustakides, ``Optimal stopping times for detecting changes in
  distributions,'' {\em Ann. Statist.}, vol.~14, pp.~1379--1387, Dec. 1986.

\bibitem{robe-technometrics-1966}
S.~W. Roberts, ``A comparison of some control chart procedures,'' {\em
  Technometrics}, vol.~8, pp.~411--430, Aug. 1966.

\bibitem{tart-etal-thirdord-2012}
A.~G. Tartakovsky, M.~Pollak, and A.~S. Polunchenko, ``Third-order asymptotic
  optimality of the generalized {Shiryaev-Roberts} changepoint detection
  procedures,'' {\em Theory of Prob and App.}, vol.~56, no.~3, pp.~457--484,
  2012.

\bibitem{lord-poll-nonantiest-2005}
G.~Lorden and M.~Pollak, ``Nonanticipating estimation applied to sequential
  analysis and changepoint detection,'' {\em Ann. Statist.}, pp.~1422--1454,
  2005.

\bibitem{gord-Poll-signrank-1994}
L.~Gordon and M.~Pollak, ``An efficient sequential nonparametric scheme for
  detecting a change of distribution,'' {\em Ann. Statist.}, pp.~763--804,
  1994.

\bibitem{li-siri-veer-isit-2014}
Y.~Li, S.~Nitinawarat, and V.~V. Veeravalli, ``Universal sequential outlier
  hypothesis testing,'' in {\em IEEE International Symposium on Information
  Theory (ISIT)}, pp.~3205--3209, 2014.

\bibitem{hero-bala-submitted-2015}
A.~Hero and B.~Rajaratnam, ``Foundational principles for large scale inference:
  Illustrations through correlation mining,''
\newblock submitted.

\bibitem{lai-ieeetit-1998}
T.~L. Lai, ``Information bounds and quick detection of parameter changes in
  stochastic systems,'' {\em IEEE Trans. Inf. Theory}, vol.~44, pp.~2917
  --2929, Nov. 1998.

\end{thebibliography}

\end{document}